\documentclass[12pt,letterpaper]{amsart}
\usepackage{graphicx}
\usepackage[utf8]{inputenc}
\usepackage{amssymb}
\usepackage{epstopdf}
\usepackage{amsmath}
\usepackage{xcolor}

\def\N{{\mathbb N}}
\def\Z{{\mathbb Z}}

\def\R{{\mathbb R}}

\def\bb{\begin}
\def\bc{\begin{center}}
\def\ec{\end{center}}
\def\be{\begin{equation}}
\def\ee{\end{equation}}
\def\ba{\begin{array}}
\def\ea{\end{array}}
\def\bea{\begin{eqnarray}}
\def\eea{\end{eqnarray}}
\def\beaa{\begin{eqnarray*}}
\def\eeaa{\end{eqnarray*}}
\def\hh{\!\!\!\!}
\def\EM{\hh &   &\hh}
\def\EQ{\hh & = & \hh}

\def\LE{\hh & \le & \hh}

\def\LT{\hh & < & \hh}

\def\AND#1{\hh & #1 & \hh}

\def\la{\lambda}

\def\da{\delta}

\def\si{\sigma}

\def\Om{\Omega}
\def\ol{\overline}
\def\ul{\underline}

\def\ck{\check}

\def\d{\cdot}

\def\dd{\cdots}
\def\oo{\infty}
\def\f{\frac}
\def\pa{\partial}
\def\z{\left}
\def\y{\right}
\def\q{\quad}
\def\qq{\qquad}

\def\dto{\downarrow}
\def\tm{\times}
\def\rd{\,{\rm d}}
\def\dt{\,{\rm d}t}

\def\dv{\,{\rm d}v}
\def\dx{\,{\rm d}x}

\def\bfl{{\bf L}}
\def\bfm{{\bf M}}

\def\bfw{{\bf W}}
\def\bfy{{\bf Y}}
\def\bfz{{\bf Z}}

\def\mcc{{\mathcal C}}

\def\mcl{{\mathcal L}}
\def\mcm{{\mathcal M}}

\def\mct{{\mathcal T}}

\def\Mo{\mcm_0}
\def\Lo{\mcl^1}

\def\Li{\mcl^\infty}

\def\bor{B_{r}}
\def\sor{S_{r}}

\def\mimr{{\bf m}_{i,m}(r)}

\def\imr{{i,m,r}}
\def\rmi{{r,m,i}}
\def\im{{i,m}}

\def\nlo{\|\cdot\|_1}

\def\nmm#1{\|#1\|}
\def\nmo#1{\|#1\|_1}

\def\nmi#1{\|#1\|_\infty}
\def\nmv#1{\|#1\|_{\bf V}}

\def\nn{\nonumber}

\def\iff{\ \Longleftrightarrow\ }

\def\andq{\quad \mbox{ and } \quad}

\def\Fr{Fr\'echet\ }

\def\ifl{\iffalse}
\def\Proof{\noindent{\bf Proof} \quad}
\def\qed{\hfill $\Box$ \smallskip}

\def\tpi{{\tilde{\pi}}}
\def\ta{\tilde{\tau}}
\def\T{\tilde{\mathcal{T}}}

\def\piim{{R}_{i,m}}
\def\hr{\hat{r}}
\def\kr{\check{r}}

\def\bb{\begin}
\def\bc{\begin{center}}
\def\ec{\end{center}}
\def\be{\begin{equation}}
\def\ee{\end{equation}}
\def\ba{\begin{array}}
\def\ea{\end{array}}
\def\bea{\begin{eqnarray}}
\def\eea{\end{eqnarray}}
\def\beaa{\begin{eqnarray*}}
\def\eeaa{\end{eqnarray*}}
\def\hh{\!\!\!\!}
\def\EM{\hh &  &\hh}
\def\EQ{\hh & = & \hh}

\def\LE{\hh & \le & \hh}

\def\LT{\hh & < & \hh}

\def\x#1{(\ref{#1})}
\def\R{{\mathbb R}}

\def\N{{\mathbb N}}

\def\nn{\nonumber}
\def\oo{\infty}
\def\lb{\label}

\begin{document}

\title[]{Explicit sharp bounds for all nodes of Sturm-Liouville operators with potentials in $L^1$ balls}

\author[]
{Jifeng Chu$^{1}$, \quad Shuyuan Guo$^2$,\quad Gang Meng$^3$,\quad Meirong Zhang$^2$}

\address{$^1$ School of Mathematics, Hangzhou Normal University, Hangzhou 311121, China}
\address{$^2$ Department of Mathematical Sciences, Tsinghua University, Beijing 100084, China}
\address{$^3$ School of Mathematical Sciences, University of Chinese Academy of Sciences, Beijing 100049, China}

\email{jifengchu@126.com (J. Chu)}
\email{guosy@mail.neea.edu.cn (S. Guo)}
\email{menggang@ucas.ac.cn (G. Meng)}
\email{zhangmr@tsinghua.edu.cn (M. Zhang)}
\thanks{Jifeng Chu was supported by the National Natural Science Foundation of China (Grant
No. 12571168), the Science and Technology Innovation Plan of Shanghai (Grant No. 23JC1403200) and Zhejiang Provincial Natural Science Foundation of China (No. Z26A010002).
Gang Meng was supported by the National Natural Science Foundation of China (Grant No. 12271509 and 12471174) and
the Fundamental Research Funds for the Central Universities (No. 5172023000540). Meirong Zhang was supported by the National Natural Science Foundation of China (Grant No. 11790273)}



\begin{abstract}
For the classical Sturm-Liouville operators,
we prove the sharp bounds for all nodes of eigenfunctions by regarding these nodes as
nonlinear functionals of potential $q\in L^1[0,1]$. By studying the optimization problems to minimize or to maximize
the nodes $\{ T_{i,m}\}$ subject to the constraint $\|q\|_{1}=r$ with $r>0$ and using the strong continuity of the nodes
in potentials, we obtain the explicit expressions for the sharp bounds, which are given as elementary functions.
\end{abstract}

\maketitle

\section{Introduction}

The study of Sturm-Liouville problems was originally started from the seminal papers \cite{s-1,s-2,sl} by Sturm and Liouville in the 1830's  and now
it has become an important and wide branch in the theory of differential equations. Since the number of the literature is vast, here we only refer the reader to
\cite{lu} for the historical account and the monographs \cite{bbw, z, ze} for detailed discussions. Let us consider the classical Sturm-Liouville
\be \lb{slp}- y'' + q(x) y =\la y,\ee
with Dirichlet boundary condition
\be \lb{bc}y(0)=y(1)=0,\ee
where $q\in \mathcal{L}^1:=L^1[0,1]$ is an integrable real potential.
It follows from the classical spectral theory \cite{bbw, z, ze} that problem \x{slp}-\x{bc} admits a sequence of real eigenvalues
$\la_m= \la_m(q)$ such that
\[\la_0<\la_1<\cdots<\la_m<\cdots,\qq \la_m\to\infty \q (m\to \infty),\]
and for each $\la_m$, there exists a corresponding real eigenfunction
\[E_m(x)=E_m(x;q), \qq m=0, 1,2,\dd.\]
Without loss of generality, we normalize the eigenfunction $E_m(x)$ such that
\[E'_m(0)>0,\qquad \int_0^1E_m^2(x)dx =1.\]
It is well-known that $E_0(x)$ has no zero in $(0,1)$ and for $m\in \N= \{1,2,\dd\}$, $E_m(x)$ has precisely $m$ zeros
\[0< T_{1,m} < T_{2,m}< \dd < T_{m,m} <1,\]
which are called the {\it nodes} of problem \x{slp}-\x{bc}.
Note that the eigenvalues $\la_m=\la_m(q)$, the eigenfunctions $E_m(x)=E_m(x;q)$ and the nodes $T_\im =T_\im(q)$ are
determined by the potential $q$ in an implicit way. Such nodes $\z\{T_\im\y\}$ are significant physical quantities. For example,
in vibration dynamics, a node represents a position in the system where no vibration occurs. See \cite{am, gzm, gzmm, st} for
qualitative properties of nodes for second order ordinary differential operators. Moreover, the information of these nodes is important in the inverse spectral theory. See for example  \cite{cmx, hmi, hwxz, mcl}.

The aim of this paper is to find the explicit sharp upper and lower bounds for all nodes of
problem \x{slp}-\x{bc} for given $L^1$ norm of the potential $q\in \Lo$.
There are already many results on sharp bounds for eigenvalues, eigenvalues ratio or eigenvalues gap for different kinds of differential operators in the literature, see for example
\cite{as, ab, cm, cmz, cmz-adv, ek, gmyz, ks, wmz, zj}. Such optimization problems have been attractive and interesting
in mathematical analysis, especially in spectral theory and the theory of differential equations. Such a topic was originally started
from the pioneering work of Krein \cite{kr}, where he studied the sharp estimates for Dirichlet eigenvalues of
\[y''+\lambda q(x)y=0\] by restricting that $q\in L^{\oo}[0,1],0\leq q\leq h,\int_0^1q(x)dx=r$ with $0<r\leq h.$

However, as far as we know, up to now there are only few results \cite{cjm, cmwz} on sharp bounds of nodes for ordinary differential operators.
For given $L^p$ norms of the potentials $q\in L^p[0,1]$  with $1<p <\oo$, it was obtained in a recent work \cite{cmwz} the optimal
characterizations of locations for all nodes of problem \x{slp}-\x{bc}. In particular, by regarding the $i$th node
of the $m$th eigenfunction as a functional of the potential, the critical equations to determine the minimizing potential
are deduced so that the node is minimized, and from the critical equations, two equivalent characterizations as nonlinear
systems for $4$-dimensional or $2$-dimensional parameters have been proved. Theoretically speaking, by taking the limit $p\downarrow 1$ of the results in \cite{cmwz}, we can obtain the results on sharp bounds for the nodes when the $L^1$ norms of the potentials are given. However,
it is very difficult because the explicit expressions of the sharp bounds cannot be obtained in \cite{cmwz} and the results in \cite{cmwz} contain several singular integrals. Besides,  we are not able to obtain the critical equations when the $L^1$-norms of the potentials are fixed, because the space $\mathcal{L}^1$ lacks local compactness even in the weak topology. For these reasons, in this paper we aim to develop a different method to work out the explicit expressions for the sharp bounds of all nodes.

Now we state our optimization problems.  For $r\in(0,\oo)$, let
$\bor:=\z\{q\in \Lo: \nmo q \le r \y\}$ and $\sor:=\z\{q\in \Lo: \nmo q = r \y\}$
be the ball and the sphere in $\Lo$ with center $0$ and radius $r$, respectively. For $m\in \N$ and $1\le i\le m$, we consider
the following optimization problems on nodes
\be \lb{tmin}\mimr :=\inf_{q\in \bor} T_\im(q) \equiv \inf_{q\in \sor} T_\im(q),\ee
\be \lb{n-max}{\bf n}_\im(r) :=\sup_{q\in \bor} T_\im(q) \equiv \sup_{q\in \sor} T_\im(q).\ee
The  reflection $x\mapsto 1-x$ leads to
\[ {\bf n}_\im(r)=\sup_{q\in \bor} T_\im(q)=1-\inf_{q\in \bor} T_{m-i+1,m}(q)=1-{\bf m}_{m-i+1,m}(r).\]
It is obvious that the solutions of \x{tmin} and \x{n-max} can lead to the following optimal estimates on the locations of nodes
\[{\bf m}_\im(\nmo q) \le T_\im(q) \le 1-{\bf m}_{m-i+1,m}(\nmo q),\quad q\in \Lo.\]
However, it is difficult to solve the optimization problem \x{tmin}.  To achieve our aim and to obtain the sharp bounds for all nodes of problem \x{slp}-\x{bc}, we need three known results on the eigenvalues
$\{\lambda_m\}$ and related eigenfunctions $\{E_m(x)\}$. The first known result is related to optimal estimates of eigenvalues $\{\lambda_m\}$,
proved by Wei, Meng and Zhang in \cite[Theorem 1.3 and Theorem 1.5]{wmz}, which will be stated as Lemma 2.5.
The second known result is about a strong continuity of $\{\lambda_m\}$ and the nodes $\z\{T_\im\y\}$ related to the potentials. That
is,  the functionals
\[(\mathcal{L}^1,w_1)\rightarrow\R:\quad q\rightarrow\lambda_m(q),\quad q\rightarrow T_\im(q)\]
are continuous, where $w_1$ is the topology of weak convergence in the space $(\mathcal{L}^1,\|\cdot\|_1)$.
Thirdly, some basic facts of measure differential equations are needed because the solution of the related optimization problems
\x{tmin} and \x{n-max} will lead to more general distributions of potentials
which have no densities with respect to the Lebesgue measure, and thus it is a very natural way to
understand such problems in the setting of measure differential equations. All of the above three known results will be presented as preliminary
results in Section 2. In Section 3, we first trickily reduce the infinite-dimensional minimization problem \x{tmin} to a finite-dimensional
optimization problem \x{Timr}. Then we obtain the solutions of the finite-dimensional minimization problem \x{Timr} by delicate analysis, and finally
we show that the solutions for the minimization problem \x{Timr} are just those for our minimization problem \x{tmin}, based on the
strong continuity of the nodes $\z\{T_\im\y\}$ related to the potentials.

\section{Preliminary results}
\setcounter{equation}{0} \lb{node-ev}

In this section, we present some preliminary known results on eigenvalues and nodes, including
sharp bounds on eigenvalues, scaling equalities for nodes and some facts on measure
differential equations.

\subsection{Measure differential equations}
We need some results on eigenvalues and nodes of linear measure differential equations \cite{m, mz, zwm}. As usual, we use $\Mo=\Mo(I):= (C(I),\nmi\d)^*$ to denote the space of (normalized Radon) measures on $I = [0,1]$. Besides the usual topology induced by the norm $\nmv{\d}$ of total variations, $\Mo$ also has the weak$^*$ topology $w^*$ which is defined by using the weak$^*$ convergence for measures. In the $w^*$ topology, the following conclusion holds.
\bb{lem} \lb{wstarlemma}   \cite{m} Let $\mu_0 \in \Mo.$ Then there exists a sequence of smooth measures $\{\mu_n \} \subset
\Mo \cap \mcc^\oo$ such that $ \|{\mu_n}\|_{\bf V} = \|{\mu_0}\|_{\bf V}  $ and
    $$ \mu_n \to \mu_0\mbox{ in $( \Mo,w^*)$}.$$
\end{lem}

Given a measure $\mu\in \Mo$. We consider the eigenvalue problem of the corresponding measure differential equations
\be \lb{ev-mu}\z\{\ba{ll}
-\rd y^\bullet(x) + y(x) \rd\mu(x)=\la y(x) \dx,\\
y(0)=0, \quad y(1)=0.
\ea \y.\ee
The structure of eigenvalues of problem \x{ev-mu} is similar to that of ordinary differential equations, i.e. problem \x{ev-mu} admits a real sequence of eigenvalues
\[\la_0(\mu) < \la_1(\mu) <\dd < \la_m(\mu) < \dd, \qq \la_m(\mu) \to +\oo \q (m\to \oo).\]
See \cite{mz, zwm} for detailed discussions. Associated with each $\la_m(\mu)$, let $E_m(x;\mu)$ be a normalized eigenfunction of problem \x{ev-mu} such that $E_m(x;\mu)$ has the $L^2$ norm $1$ and is positive for $x$ in some pinched neighborhood $(0,\da)$. For each $m\in \N$, $E_m(x;\mu)$ admits $m$ nodes
\[0< T_{1,m}(\mu) < T_{2,m}(\mu) < \dd < T_{m,m}(\mu) <1.\]
Moreover, we have the following strong continuity by regarding such nodes as functionals in the measure $\mu\in \Mo$.

\bb{lem} \lb{CC} {\rm\cite{gzm, gzmm}}
Let $m\in \N$ and $1\le i\le m$ be given. By regarding the node $T_\im(\mu)$ of problem {\rm\x{ev-mu}} as a functional from $\Mo$ to $(0,1)\subset\R$, it admits a strong continuity in the following sense: $T_\im(\mu_n) \to T_\im(\mu)$ if
$\mu_n \to \mu$  in $(\Mo,w^*)$.  In particular, the node $T_\im(q)$ for problem {\rm\x{slp}-\x{bc}} are  continuous in potential $q\in (\mathcal{L}^1,w_1)$.

\end{lem}

\bb{lem}\lb{CD} {\rm \cite{gxz}}
Let $m\in \N$ and $1\le i\le m$ be given. Then $T_\im(q)$ is continuously \Fr differentiable in $q\in(\Lo,\nlo)$.
Moreover, the \Fr derivative is given as
\be \lb{Fr}\pa_q T_\im(q)(x)= H_\im(x;q) E_m^2(x;q) \in \Li,\ee
where $H_\im(x;q)$ is the two-step function
\[H_\im(x;q)= \z\{ \ba{ll} +a_\im, &  x\in[0,T_\im(q)), \\
- b_\im, & x\in[T_\im(q),1],\ea \y.\]
and the two constants $a_\im=a_\im(q),~b_\im=b_\im(q)$ are given as
\[a_\im = \int_{[T_\im(q),1]} \z(\f{E_m(v;q)}{E'_m(T_\im(q);q)}\y)^2\dv, \q
b_\im = \int_{[0,T_\im(q)]} \z(\f{E_m(v;q)}{E'_m(T_\im(q);q)}\y)^2\dv.\]
\end{lem}

The following result on $\mimr$ can be easily obtained from Lemma \ref{CC} and Lemma \ref{CD}.

    \bb{lem} \lb{TerP1}
For $r\in(0,\oo)$, $m\in\N$ and $1\le i\le m$, it holds that    \[
    \mimr \in \z(0, i/(m+1)\y).
    \]
    \end{lem}

\Proof For the zero potential, one has $T_\im(0) = i/(m+1)$. Let us choose a perturbed potential
    \[
    q_t(x) =\z\{ \ba{ll} - t, & \mbox{for } x\in[0,i/(m+1)], \\
    0 , & \mbox{for } x\in(i/(m+1),1],
    \ea
    \y.
    \]
where $0<t \ll 1$. Then $q_t \in B_r$ and it follows from the formula \x{Fr} that
    \[
    \f{\rd}{\dt} T_\im(q_t) = -\int_{[0,i/(m+1)] } a_\im E_m^2(x;0) <0.
    \]
Thus $T_\im(q_t) < i/(m+1)$ for $0<t \ll 1$. Therefore $\mimr < i/(m+1)$.

On the other hand, let $\{q_n\}\subset S_r$ be a minimizing sequence such that
    \[
    \lim_{n\to\oo} T_\im(q_n) = \mimr.
    \]
Let $\{q_n\}$ be the densities of the absolutely continuous measures
    \be \lb{muq}
    \mu_{q_n}(x) = \int_{[0,x]} q_n(t) \dt, \q x\in I.
    \ee
Then $\nmv{\mu_{q_n}} =r$. By the Banach theorem, without loss of generality, we can assume that $\{\mu_{q_n}\}$ has a weak$^*$ limit measure $\mu_r$. Thus it follows from Lemma \ref{CC} that
    \[
    \mimr = \lim_{n\to\oo} T_\im(q_n) =\lim_{n\to\oo} T_\im(\mu_{q_n}) = T_\im(\mu_r).
    \]
In particular, $\mimr = T_\im(\mu_r)>0$.
\qed

\subsection{Known results on eigenvalues and nodes}
Let $r\in(0,\infty)$ and $m\in\Z_+$ be given. Consider the following optimization  problems related to the eigenvalues
\be \lb{lapr} \bfl_{m}(r):=\inf_{q\in S_{r}} \la_{m}(q) = \inf_{q\in B_{r}} \la_{m}(q),\quad
\bfm_{m}(r):= \sup_{q\in S_{r}} \la_{m}(q) = \sup_{q\in B_{r}} \la_m(q).\ee

\bb{lem}\lb{N50} {\rm \cite{wmz}}
Let $r\in(0,\infty)$ be given. Then
\[\bfl_{0}(r) = \la_0\z(-r \da_{1/2}\y)=\bfz^{-1}(r),\qq
\bfm_{0}(r)= \la_0\z(\eta_{r}\y)=\bfy(r),\]
where $\bfz$ is a decreasing diffeomorphism from $\left(-\infty, \pi^{2}\right]$ to $[0, +\infty)$ defined by
\be \lb{bfz}\bfz(x)=\left\{\begin{array}{ll}
2 \sqrt{-x} \operatorname{coth} (\sqrt{-x} / 2) & \text { for } x \in(-\infty, 0), \\
4 & \text { for } x=0, \\
2 \sqrt{x} \cot (\sqrt{x} / 2) & \text { for } x \in (0, \pi^{2}],
\end{array}\right.\ee
and $\bfy:\left[0,+\infty\y) \rightarrow\z[\pi^2, +\infty\y)$ is an increasing diffeomorphism given as
\be \lb{bfy}\bfy(x)=\z(\f{\pi + \sqrt{\pi^2 + 4x}}{2}\y)^2,    \ee
where $\kappa \da_a$ denotes the Dirac measure of mass $\kappa$ at point $a$ and  $\eta_r$ is the potential defined as
\[\eta_r(x):=\left\{\begin{array}{ll}\mathbf{Y}(r),
&  x\in[\tau_r,1-\tau_r],  \\
0, &  x\in [0, \tau_r) \cup (1-\tau_r,1],
\end{array}\right.\]
with the constant
\[\tau_{r}=\frac{\pi}{2 \sqrt{\mathbf{Y}(r)}} \in\left(0,1/2\right).\]
Moreover, for any $m\in \N$, we have the following the dilation relations
\be \lb{dr-1}\bfl_{m}(r) = (m+1)^2\bfl_{0}\z(\f{r}{(m+1)^2}\y)= (m+1)^2\bfz^{-1}\z(\f{r}{(m+1)^2}\y) = \la_m(\mu_{m, r}),\ee
\be \lb{dr-2}\bfm_{m}(r) = (m+1)^2\bfm_{0}\z(\f{r}{(m+1)^2}\y)= (m+1)^2\bfy\z(\f{r}{(m+1)^2}\y) = \la_m(\eta_{m, r}),\ee
where \be \lb{mumr}
 \mu_{m, r }(x) = -\frac{r}{m+1} \sum_{k=0}^{m} \delta_{\frac{2k+1}{2(m+1)}}(x),
    \ee
and
\be \lb{etamr}
    \eta_{m,r}(x):=\left\{\begin{array}{ll}
   (m+1)^2\mathbf{Y}({r}/{(m+1)^2}),
    &  x\in[\frac{\tau_{m,r}+k}{m+1},\frac{1-\tau_{m,r}+k}{m+1}],\ k = {0,1,\dd,m},    \\
    0, & \mbox { elsewhere},
    \end{array}\right.
    \ee
with
\[\tau_{m,r}\equiv\frac{\pi}{2 \sqrt{\mathbf{Y}({r}/{(m+1)^2})}} \in\left(0, \f{1}{2}\right).\]

\end{lem}

The above conclusions  can be extended from the interval $[0,1]$ to any interval $[x_1,x_2]$ with $x_2-x_1 = T>0.$
Indeed, suppose  $r >0$ and
\[ B_{r,T}:=\z\{q\in \Lo(x_1,x_2): \|q\|_{1,[x_1,x_2]}  \le r \y\}.\]
 Let   $\la_m(q,T)$ be the   lowest eigenvalues  of equation \x{slp}   with Dirichlet boundary condition
 \[y(x_1) = y(x_2) =0.\]
 By the similar arguments, we can  solve the corresponding optimization problems as follows
\be \lb{L1rBT} \bfl_m(r,T):= \inf \Big\{ \la_m(q,T): \mu \in B_{r,T}\Big\} = \la_m\Big(\mu_{m, r }( \f{x-x_1}{x_2-x_1})\Big)= \f{(m+1)^2}{T^2}\bfz^{-1}\z( \f{T r}{(m+1)^2} \y ),\ee
\be \lb{M1rBT}\bfm_m(r,T):= \sup \Big\{ \la_m(q,T): \mu \in B_{r,T}\Big\} = \la_m\Big(\eta_{m,r}( \f{x-x_1}{x_2-x_1})\Big)= \f{(m+1)^2}{T^2}\bfy \z (\f{T r}{(m+1)^2} \y).\ee

\subsection{Scaling equalities for nodes and eigenvalues}

Let $m\in \N$ and $1\le i\le m$ be given. For a potential $q\in\Lo$, we have the node $T=T_\im(q)\in (0,1)$. By using $T$ to cut the interval $[0,1]$ into two subintervals $[0,T]$ and $[T,1]$, we obtain the following two Dirichlet eigenvalue problems on subintervals
\be \lb{ev0t}-E_m^{\prime \prime}(x)+q(x)E_m(x) = \lambda_{m}(q)E_m(x), \qq E_m(0)=E_m(T)=0,\ee
\be \lb{evt1}-E_m^{\prime \prime}(x)+q(x)E_m(x) = \lambda_{m}(q)E_m(x), \qq E_m(T)=E_m(1)=0.\ee
By using the scaling $x\mapsto T x$ and the scaling $x\mapsto T+(1-T)x$ for \x{ev0t} and \x{evt1} respectively, we obtain the following two Dirichlet eigenvalue problems
\be \lb{haty}-\hat{y}^{\prime \prime}(x) +\hat{q}(x) \hat{y}(x)= T^{2} \lambda_{m}(q)\hat{y}(x), \qq
\hat{y}(0)=\hat{y}(1)=0,\ee
\be \lb{cky}-\ck{y}^{\prime \prime}(x) +\ck{q}(x) \ck{y}(x)= (1-T)^{2} \lambda_{m}(q) \ck{y}(x), \qq
\ck{y}(0)=\ck{y}(1)=0,\ee
where
\be \lb{qq}\hat{q}(x):=T^{2} q(Tx)\in\Lo,\qq \ck{q}(x):=(1-T)^{2} q(T+(1-T)x)\in\Lo,\ee
\[\hat{y}(x)=\hat{y}(x;q):=E_m(Tx),\qq \ck{y}(x)=\ck{y}(x;q):=E_m\left(T+(1-T)x\right).\]
Since $\hat{y}(x)$ has exactly $i-1$ interior zeros in $[0,1]$ and $\ck{y}(x)$ has exactly $m-i$ interior zeros, we can deduce from \x{haty} and \x{cky} that $\hat{q}$ and $\ck{q}$ have Dirichlet eigenvalues
\be \lb{qla12}\lambda_{i-1}\left(\hat{q}\right)=T^{2} \lambda_{m}(q),\qq\lambda_{m-i}\left(\ck{q}\right)=(1-T)^{2} \lambda_{m}(q),\ee
which can  produce the following scaling equality between nodes and eigenvalues
\be\lb{la1la2}\frac{\la_{i-1}(\hat{q})}{T^2} = \frac{\la_{m-i}\z(\ck{q}\y)}{(1-T)^2}.\ee
Note that the two equalities in \x{qla12} are the usual scaling equalities for Dirichlet eigenvalues.

It is easy to see that the converse of the above construction is also true. Thus we have obtained the following result.

\bb{lem}\lb{revert}
Let $m\in \N$ and $1\le i\le m$ be given. Then for any $q\in \Lo$, one has
\[\z(T,\hat{q},\ck{q}\y)= (T_\im(q), \hat{q}(\d;q), \ck{q}(\d;q)),\]
which satisfies the scaling equality {\rm\x{la1la2}.} Conversely, let $\hat{q}$, $\ck{q} \in\Lo$ and $T\in(0,1)$ satisfy the scaling equality {\rm\x{la1la2}.} Then there must be a potential $q\in\Lo$ such that
$(T_\im(q), \hat{q}(\d;q), \ck{q}(\d;q))=\z(T,\hat{q},\ck{q}\y).$
\end{lem}

\section{Main results}
\setcounter{equation}{0} \lb{Inp1}

In this section, we deduce the infinite-dimensional minimization problem \x{tmin} to a finite-dimensional optimization problem. Then
we prove the solutions of the finite-dimensional problem and finally show that such solutions are just those of our infinite-dimensional minimization problem \x{tmin}.

\subsection{A related finite-dimensional optimization problem} For given $r\in(0,\oo),~m\in\N, ~1\le i\le m,$ we consider any integrable potential $q\in \sor$ such that the node $T=T_\im(q)\in(0,i/(m+1))$.
Let us introduce two non-negative parameters
\[\hr=\hr(q):=\nmm{q}_{\Lo{[0,T]}},\quad \kr=\kr(q):=\nmm{q}_{\Lo{[T,1]}}.\]
It is trivial that
\be \lb{r1r2}\hr + \kr =r.\ee
Moreover, by the equalities \x{qq} for scaling potentials $\hat{q}$ and $\ck{q}$, one has
\[\nmo{\hat{q}} = T \hr, \qq \nmo{\ck{q}} = (1-T) \kr.\]
From the optimization problems \x{lapr} on eigenvalues, one has
\[\la_{i-1}(\hat{q})\ge \bfl_{i-1}(\nmo{\hat q})=\bfl_{i-1}(T \hr), \quad\la_{m-i}(\ck q) \le \bfm_{m-i}(\nmo{\ck q})= \bfm_{m-i}((1-T) \kr).\]
It follows from the scaling result \x{la1la2} that we have the following inequality
\be\lb{ine}\f{\bfl_{i-1}(T \hr)}{T^2} \le  \f{\bfm_{m-i}\z((1-T) \kr\y)}{(1-T)^2},\ee
which contains three quantities $(\hr, \kr, T)$. Motivated by \x{r1r2} and \x{ine}, we introduce the following function
\be \lb{ftrfun}f(T,\hr):= \bfl_{i-1}\z(T \hr\y)- \z(\frac{T}{1-T}\y)^2\bfm_{m-i}\Big((1-T)(r-\hr)\Big),\q (T,\hr)\in {\mathcal D}:=(0,i/(m+1))\tm [0,r].\ee
Next we focus on the problem of the implicit function $f(T,\hr)=0$, which is equivalent to solve the equation
\be\lb{P2}\bfl_{i-1}\z(T \hr\y) = \z(\frac{T}{1-T}\y)^2\bfm_{m-i}\Big((1-T)(r-\hr)\Big).\ee
Since both $\bfl_{i-1}$ and $\bfm_{m-i}$ are explicitly given in Lemma \ref{N50}, equation \x{P2} can be used to determine $T$ from $\hr$ in an implicit way.

\bb{lem} Given $r\in(0,\oo),~m\in\N, ~1\le i\le m.$ Then there exists a unique $C^1$ function
\be \lb{PTr1}T=\mct(\hr) =\mct_\imr(\hr):[0,r]\to (0,i/(m+1)),\ee
such that
\be \lb{PTr12}\bfl_{i-1}\z(\mct(\hr)\hr\y) = \z(\frac{\mct(\hr)}{1-\mct(\hr)}\y)^2\bfm_{m-i}\Big((1-\mct(\hr)) (r-\hr)\Big),\quad \hr\in[0,r].\ee
\end{lem}

\Proof For any fixed $\hr\in[0,r]$, one has
\[\lim_{T\dto 0} f(T,\hr) = \bfl_{i-1}\z(0\y) = (i\pi)^2 >0,\]
and \beaa
f\z(\f{i}{m+1},\hr\y) \EQ \bfl_{i-1}\z( \f{i\hr}{m+1}\y) -\f{\z(\f{i}{m+1}\y)^2}{\z(1- \f{i}{m+1}\y)^2} \bfm_{m-i}\z( \f{i(r-\hr)}{m+1}\y) \\
\EQ (i\pi)^2 \z( \f{\bfl_{i-1}\z( \f{i\hr}{m+1}\y)}{(i\pi)^2} - \f{\bfm_{m-i}\z(\f{i(r-\hr) }{m+1}\y)}{((m-i+1)\pi)^2} \y) .
\eeaa
By the definition of the optimization problems \x{lapr}, one has
\[\f{\bfl_{i-1}\z(\f{i\hr }{m+1}\y)}{(i\pi)^2}\le 1, \qq  \f{\bfm_{m-i}\z(\f{i(r-\hr) }{m+1}\y)}{((m-i+1)\pi)^2}\ge 1.\]
Note that at least one of the inequalities above is strict because at least one of $\{\hr,\kr\}$ is not zero. Hence
\[f(i/(m+1),\hr)<0.\]
Therefore there exists at least one solution $T=\mct(\hr)=\mct_\imr(\hr)\in (0,i/(m+1))$ such that
\[ f(\mct(\hr),\hr) =0.\]

To prove the uniqueness, now we show that
\be \lb{dft}\f{\pa }{\pa T}f(T,\hr)<0.\ee
In fact, it follows from the dilation relations in \x{bfz} and \x{dr-1}-\x{dr-2} that
\beaa\bfl_{i-1}(T\hr) \EQ i^2 \bfz^{-1}\z({T\hr}/{i^2}\y),\\
\f{\pa}{\pa T}\bfl_{i-1}(T\hr) \EQ \hr\d \z(\bfz^{-1}\y)'\z({T\hr}/{i^2}\y)\le 0,\eeaa
since the function $\bfz$ is decreasing. Moreover,
\beaa\EM \bfm_{m-i}\Big((1-T)(r-\hr)\Big) = (m-i+1)^2\bfy \z(\f{(1-T)(r-\hr)}{(m-i+1)^2} \y)\\
\EQ \z(\f{(m-i+1)\pi + \sqrt{((m-i+1)\pi)^2 + 4 (r-\hr)(1 -T)}}{2}\y)^2.
\eeaa
Consider
\beaa
g(T,\hr) \AND{:=} \sqrt{\z(\frac{T}{1-T}\y)^2\bfm_{m-i}\Big((1-T)(r-\hr)\Big)}\\
\EQ \f{1}{2} \frac{T}{1-T}\z((m-i+1)\pi + \sqrt{((m-i+1)\pi)^2 + 4 (r-\hr)(1 -T)}\y).
\eeaa
Then for all $(T,\hr)\in {\mathcal D}$,
\[\f{\pa }{\pa T}g(T,\hr) = \frac{1}{2(1-T)^2} \z((m-i+1)\pi+
\frac{((m-i+1)\pi)^2+2(r-\hr)(1-T)(2-T)}{\sqrt{((m-i+1)\pi)^2 + 4(r-\hr)(1 -T)}} \y) > 0.\]
As a consequence,
\[\f{\pa}{\pa T}f(T,\hr) = \f{\pa}{\pa T}\bfl_{i-1}(T\hr)- 2g(T,\hr) \d \f{\pa }{\pa T}g(T,\hr)<0.\]
Thus \x{dft} has been verified. In particular, when $\hr\in [0,r]$ is fixed, $f(T,\hr)$ is strictly decreasing in $T\in(0,1)$. Hence
we have proved the uniqueness of the implicit function $\mct_\imr(\hr)$ satisfying \x{PTr12}.
Finally, due to \x{dft} and by the implicit function theorem, $\mct_\imr(\hr)$ is smooth in $\hr$.  \qed

Let us analyze the dependence of the function $\mct_\imr(\hr)$ on $\hr$ as well as on $(\imr)$, including the monotonicity of $\mct_\imr(\hr)$ in $\hr$. Since $\mct_\imr(\hr)$ is defined in an implicit way, we introduce a new independent variable $t$ for replacing $\hr$. For $\hr\in[0,r]$ and $\kr=r-\hr\in[0,r]$, we know that $T= \mct(\hr) = \mct_\imr(\hr)\in(0,i/(m+1))$. By using the relations \x{dr-1} and \x{dr-2}, the equality \x{P2} for defining $T=\mct_\imr(\hr)$ is
\be \lb{Lg0}
\bfz^{-1}\left(\f{\mct(\hr)\cdot \hr}{i^2}\right)= \z( \f{\mct(\hr)}{1-\mct(\hr)} \y)^2 \z( \f{m-i+1} {i} \y)^2\bfy\z( \f{(1-\mct(\hr)) \d \kr}{(m-i+1)^2}\y)>0.\ee
It follows from the formulas \x{bfz} that
\[\f{\mct(\hr)\cdot \hr}{i^2} \in [0,4),\quad \f{\bfl_{i-1}\left(\mct(\hr)\cdot \hr\right)}{ i^2}\in (0,\pi^2].\]
Note that $\bfz^{-1}(\la)$ is decreasing in $\la\in (-\oo,\pi^2]$, we can introduce the parameter
\be \lb{tr1}t:=t\left(\hr\right)=\frac{\sqrt{\bfl_{i-1}\left(\mct(\hr)\cdot \hr\right)}}{2i} =\f{1}{2} \sqrt{\bfz^{-1}\left(\f{\mct(\hr)\cdot \hr}{i^2}\right)} \in (0,\pi/2].\ee
Now we have three parameters $\hr$, $t$ and $\mct\ (=\mct(\hr))$, which are correlated with each other by
the following way. It follows from \x{Lg0}, \x{tr1} and \x{bfz} that
\be \lb{Eq1}\mct \hr=4 i^2 t \cot t.\ee
By using \x{Lg0} and the expression \x{bfy} for $\bfy$, we obtain the second equality
\bea \lb{Eq2}t \EQ \frac{\mct}{2(1-\mct)}\d\frac{m-i+1}{2i} \cdot \left(\pi+\sqrt{\pi^{2}+ \frac{4(1-\mct)(r-\hr)}{(m-i+1)^{2}}}\right)\nn\\
\EQ \frac{\mct}{2(1-\mct)} \left(\tpi +\sqrt{\tpi ^{2}+ (1-\mct)\f{r-\hr}{i^2}}\right),
\eea
in which the constant $\tpi$ is defined as
\be \lb{pim}\tpi=\tpi_\im:= \f{m-i+1}{2i} \pi \in \z[ \f{\pi}{2m},\f{m\pi}{2}\y].\ee

Equations \x{Eq1} and \x{Eq2} for three parameters $(\hr, t, \mct)$ can be solved as follows. First, equation \x{Eq2} can be transformed to the quadratic equation for $\mct$
\[\f{r-\hr}{i^2} \mct^2 +4 t(t+\tpi) \mct -4 t^2=0,\]
which has the positive solution
\be \lb{trr1}\mct=\frac{2 t}{t+\tpi +\sqrt{\z(t+\tpi \y)^{2}+\frac{r-\hr}{i^2}}}.\ee
 Substituting \x{trr1} into \x{Eq1}, we obtain
\[ \frac{2 t\hr}{t+\tpi +\sqrt{\z(t+\tpi \y)^{2}+\frac{r-\hr}{i^2}}} =4 i^2 t \cot t, \]
which is equivalent to a quadratic equation for $\hr$
\be \lb{rr1a}\hr^{2}+ 4i^2\cot t \z(\cot t -t-\tpi \y)\hr- 4i^2 r \cot ^{2} t=0.\ee
As $\hr$ is non-negative, equation \x{rr1a} is solved as
\be \lb{r1t}\hr=\hr(t)=2i^2 \cot t \cdot\z(t-\cot t+\tpi +\sqrt{\left(t-\cot t +\tpi \right)^{2}+\frac{r}{i^2}}\y),\ee
which is an explicit expression for the inverse of $t=t(\hr)$. Finally, by using \x{r1t}, one has from \x{Eq1} that
\be \lb{trr2}\mct=\frac{2i^2 t}{r} \cdot\left(\sqrt{\left(t-\cot t+\tpi \right)^{2}+\frac{r}{i^2}} -\left(t-\cot t+\tpi \right)\right).\ee

Note that the functions $\hr=\hr(t)$ and $\mct=\mct(t)$, given as \x{r1t} and \x{trr2} above, are well-defined for $t$ in the interval $(0,\pi/2]$. Now we
compute the values of $t(\hr)$ and $\mct_\imr(\hr)$ at end-points $\hr=0$ and $\hr=r$.

\bb{lem}\lb{ine1} The values of $t(\hr)$ and $\mct_\imr(\hr)$ at end-points $\hr=0,\hr=r$ can be computed as follows
\be \lb{thr0}t(\hr)|_{\hr=0}=\pi/2,,\quad \mct(\hr)|_{\hr=0}=\f{2i \pi}{(m+1)\pi+\sqrt{((m+1)\pi)^2 + 4 r}},\ee
\be\lb{Thrr}\mct(\hr)|_{\hr=r}=\T = \T_\imr := \bfw^{-1}_\im(r)\in (0,i/(m+1)),\ee
and
\be\lb{thrr}t(\hr)|_{\hr=r}=\ta=\ta_\imr := \frac{\tpi\d\bfw_\im^{-1}(r)}{1-\bfw_\im^{-1}(r)} \in \z(0,{\pi}/{2}\y),\ee
where $\bfw_\im:(0,i/(m+1))\to (0,+\oo)$ is given by
\be \lb{WimT}\bfw_\im(T):=\f{4i^2\tpi}{1-T} \cot\f{\tpi T}{1-T}\equiv \f{2i(m-i+1)\pi}{1-T} \cot\f{(m-i+1)\pi T}{2i(1-T)},\ee
and the constants $\T, \ \ta$ in {\rm\x{Thrr}, \x{thrr}} are correlated with
\be \lb{trtr}\ta= \f{\tpi \T }{1-\T}, \quad\T=\frac{\ta}{\ta+\tpi }.\ee
Moreover, by introducing a decreasing diffeomorphism $P:  (0,\pi/2]\to [0, +\oo)$ as
\be \lb{Pim}P(t)=P_\im(t):=4 i^{2}\left(t+\tpi \right) \cot t,\ee
the constants $\tpi=\tpi_\im$ in {\rm\x{pim}} and $\ta=\ta_\imr$ in {\rm\x{thrr}} are correlated with $P(\ta)=r$, i.e.
\be \lb{corr}4 i^2 (\ta +\tpi) \cot \ta =r.\ee
\end{lem}

\Proof At $\hr=0$, since $t\in(0,\pi/2]$, we obtain from \x{Eq1} that $t=\pi/2$ and thus the first equality of \x{thr0} holds. By using \x{trr2} and \x{thr0}, the second equality of \x{thr0} is obivous.

If $\hr=r$, \x{Eq1} and \x{Eq2} can be simplified into
\be \lb{Eq3}r \mct = 4 i^2 t\cot t, \qq t= \f{\tpi \mct }{1-\mct}.\ee
By the second equation of \x{Eq3}, we have
\[\mct = \f{t}{t+\tpi},\]
and then the first equation of \x{Eq3} becomes
\[\bfw_\im(\mct) = r.\]
Hence $\mct=\bfw_\im^{-1}(r)$. Going to the second equation of \x{Eq3}, we obtain $t= \ta$. The relations in \x{trtr} are obvious from the second equation of \x{Eq3}. Note that we can rewrite \x{thrr} as
\[\bfw_\im^{-1}(r)= \f{\ta}{\ta+\tpi},\quad\bfw_\im\z( \f{\ta}{\ta+\tpi}\y) =r.\]
Then the equality \x{corr} holds by using \x{WimT} and \x{Pim}.\qed

\bb{lem}\lb{mon} The function $\hr= \hr(t)$ is a decreasing diffeomorphism from $[\ta,\pi/2]$ onto $[0,r]$. In fact, one has
\be \lb{dhrdt}\f{\rd \hr(t)}{\dt} <0,\quad t\in [\ta,\pi/2].\ee
\end{lem}

\Proof  From \x{Eq1}, $\hr=\hr(t)$ satisfies
\[\mct(\hr(t)) \hr(t)=4 i^2 t \cot t,\qq t\in[\ta,\pi/2],\]
where $\mct(\hr) = \mct_\imr(\hr)$ is a $C^1$ function. Differentiating the above equality with respect to $t$, we obtain that
\[\z(\mct(\hr(t)) + \z.\f{\rd \mct(\hr)}{\rd \hr}\y|_{\hr =\hr(t)}  \y) \d\f{\rd \hr(t)}{\dt} = 2i^2 \f{\sin(2t) - 2t}{\sin^2 t}<0,\quad t\in [\ta,\pi/2].\]
Hence $\f{\rd \hr(t)}{\dt} $ does not change sign. Using \x{r1t}, it is easy to obtain
\[\z. \f{\rd \hr(t)}{\dt}\y|_{t=\pi/2} = -  i^2 \z( {\pi}+2\tpi + \sqrt{\z({\pi}+2\tpi\y)^2+\f{4r}{i^2}}\y) <0.\]
Therefore, we have proved the desired inequality \x{dhrdt}. \qed

Recall that $\mct(t)=\mct_\imr(t)$ is given by \x{trr2}. A direct calculation gives the derivative
\be \lb{Tp1}\f{\rd \mct}{\dt}= f(t) \d h(t), \qq t\in [\ta,\pi/2],\ee
where
\[f(t)=\f{2i^2}{r} \f{\sqrt{\left(t-\cot{t}  + \tpi  \right)^{2}+\frac{r}{i^{2}}} - \left(t-\cot{t}  + \tpi  \right)}
{\sqrt{\left(t-\cot{t}  + \tpi  \right)^{2}+\frac{r}{i^{2}}}\z( \sqrt{\left(t-\cot{t}  + \tpi  \right)^{2}+\frac{r}{i^{2}}} + t \left(2+\cot^{2}{t}\right)\y)}>0,\quad t\in[\ta,\pi/2],\]
and
\be \lb{himr}h(t):=h_\imr(t)=\left(t-\cot{t}  + \tpi  \right)^{2}+\frac{r}{i^{2}}-t^{2}\left(2+\cot^{2}{t}\right)^{2},\q t\in \z[\ta, \pi/2\right].\ee

Note that the function $h$ may change sign on the interval $[\ta, \pi/2]$. Moreover, it follows from \x{Tp1} and Lemma \ref{mon} that the following result holds.

\bb{lem} \lb{signs} The functions
\[\f{\rd \mct_\imr(t)}{\dt}\andq h_\imr(t)\]
have the same sign at each $t\in[\ta,\pi/2]$, while the functions
\[\f{\rd \mct_\imr(\hr)}{\rd\hr}\andq  h_\imr(t)\]
have the reverse sign whenever $\hr$ and $t$ are corresponded via $\hr=\hr(t)$.
\end{lem}

Let us consider the following finite-dimensional optimization problem
 \be  \lb{Timr}   T_\imr:=\min_{\hr\in [0,r]} \mct_\imr(\hr) = \mct_\imr(\hr_\imr) .
  \ee
It will be shown that the solutions of the minimization problem \x{Timr} are different in two cases $i\leq \frac{m+1}{2}$
and $i>\frac{m+1}{2}$. The latter is more complicated and more precise. First we consider the case $i\leq \frac{m+1}{2}$.

\bb{lem} \lb{ile}
Assume that $(\rmi)$ satisfies
\be \lb{im1}i\le \frac{m+1}{2}.\ee
Then the function $\mct_\imr(\hr)$ is decreasing in $\hr\in [0,r]$. As a consequence, $\hr_\imr = r$ and
\be \lb{le_m}  T_\imr = \mct_\imr(r) =\T= \bfw_\im^{-1}(r).\ee
\end{lem}

\Proof Let $t\in \z[\ta, \pi/2\right)$, i.e. $\hr=\hr(t)\in(0,r]$. From \x{rr1a}, we have
\[-2\cot t \d (t+\tpi) = 2\f{r-\hr}{\hr} \cot^2 t-\f{\hr}{2i^2},\]
here $\hr=\hr(t)\in(0,r]$ is a function of $t$ and is given by \x{r1t}. It follows from \x{himr} that
\beaa
h_\imr(t)\EQ - 2 \cot{t}\d \left(t+\tpi \right) +\left(t+\tpi \right)^{2} + \cot ^{2}{t}- t^{2}\left(4+4\cot^{2}{t} +\cot^4 t\right)+\frac{r}{i^{2}}\\
\EQ 2 \frac{r-\hr}{\hr} \cot ^{2} t-\frac{\hr}{2 i^{2}}+\left(t+\tpi \right)^{2}+ \cot ^{2}{t}- t^{2}\left(4+4\cot^{2}{t} +\cot^4 t\right)+\frac{r}{i^{2}}\\
\EQ \left(\tpi -t\right)\z(\tpi+3 t \y) + \z(\frac{2r-\hr}{2 i^{2}}- 4\z(t\cot{t}\y)^2\y) +\z(1-(t\cot t)^2 + 2\f{r-\hr}{\hr} \y)\cot^2 t\\
\AND{=:} I_1 + I_2 +I_3.\eeaa
Due to condition \x{im1}, we have that
\[\tpi = \f{m-i+1}{2i}\pi \ge \f{\pi}{2}.\]
Since $t\cot t$ is decreasing in $[\ta,\pi/2]$, then
$0 < t\cot t<1$ and thus
\[I_1=\left(\tpi -t\right)\z(\tpi+3 t \y)>0,\quad t\in[\ta,\pi/2).\]
For the second term $I_2$, it follows  from \x{PTr1} and \x{Eq1} that
\[\f{\hr}{2i^2}=\frac{1}{T} \cdot 2 t \cot{t}> \f{m+1}{i} \d 2 t\cot t \ge 4 t\cot{t},\quad t\in[\ta,\pi/2).\]
As $\hr\in[0,r]$, we obtain that
\[I_2 =\frac{2r-\hr}{2 i^{2}}- 4\z(t\cot{t}\y)^2 \ge \frac{\hr}{2 i^{2}}- 4\z(t\cot{t}\y)^2 >4 t\cot{t}- 4(t\cot{t})^2 >0,\quad t\in[\ta,\pi/2).\]
In a similar way, we can also have
\[I_3=\z(1-(t\cot t)^2 + 2\f{r-\hr}{\hr} \y)\cot^2 t>0,\quad t\in[\ta,\pi/2).\]
Thus we conclude that
\[h_\imr(t)> 0,\qquad t\in[\ta,\pi/2).\]
Moreover, it can be easily checked that $h_\imr(t)=0$ if and only if
\[i=\frac{m+1}{2}, \quad t=\pi/2.\]
By Lemma \ref{signs}, we know that $\mct_\imr(\hr)$ is strictly decreasing in $\hr\in [0,r]$. Therefore $\hr_\imr=r$ and \x{le_m} holds.
\qed

Now we consider the second case
\be \lb{im2}i> \frac{m+1}{2},\ee which yields that
\[\tpi= \tpi_\im= \f{m-i+1}{2i}\pi < \f{\pi}{2}.\]
From Lemma \ref{signs}, one knows that it is important to determine the sign of the function $h_\imr$ on $[\ta,\pi/2]$.
First we consider the right end-point $t=\pi/2$. By \x{himr} and Lemma \ref{signs}, one has
\be \lb{him1} h_\imr(\pi/2)= \f{r-\piim}{i^2}, \ee
where \[\piim := \f{\pi^2}{4}\z((2i)^2-\z(m+1\y)^2\y) >0.\]

It follows from Lemma \ref{signs} and \x{him1} that the following results hold.

    \bb{lem}\lb{lem5}
Suppose that $(\im)$ satisfies {\rm\x{im2}}. Then
\beaa \lb{l51}
h_\imr\left(\pi/2\right) \leqslant 0 \iff r \leqslant \piim  \left.\iff \f{\rd \mct_\imr\left(\hr\right)}{\rd \hr}\right|_{\hr=0} \geqslant 0,\\
\lb{152}
h_\imr\left(\pi/2\right) \geqslant 0 \iff r \geqslant \piim  \left.\iff \f{\rd \mct_\imr\left(\hr\right)}{\rd \hr}\right|_{\hr=0} \leqslant 0.
\eeaa
Moreover, these equivalent inequalities take the equal sign at the same time. The symbol ``$\iff$" means ``if and only if".
\end{lem}

Besides the elementary functions $P(t)=P_\im(t)$ and $h(t)=h_\imr(t)$, we will encounter other four
elementary functions in the sequel proofs, that is,
\[U: (0,\pi/2]\rightarrow \R,\quad H: (0,\pi/2]\to \R,\quad g: (0,\pi/2]\to\R,\quad F:(0,\pi/2)\to \R\] are defined as
\be \lb{Uim}U(t): = U_\im(t) =t+t \cot ^2 t-\cot t-\tpi,\ee
\be\lb{Him}H(t)=H_\im(t):=\left(t+\cot t+\tpi \right)^{2}-\left(2t+t\cot ^{2} t\right)^{2},\ee
\be \lb{gim}g(t)=g_\im(t):= -\cot{t}+\frac{2 t^{2} \cot{t}}{\sin ^{2} t}-\frac{t}{\sin ^{2} t}+\tpi,\ee
\be \lb{Funct}F(t):=
\left(t + t \cot^{2}{t} - \cot{t} \right)
\left(t + t\cot^{2}{t} + 3\cot{t} + \pi \right) + 4t \cot{t} - \frac{3}{4}\pi^{2}.
\ee

Obviously, $U$ is an increasing function. We can also verify that $H$ is a decreasing function and $g$ is a decreasing function since
we have that by direct computations for all $ t\in(0,\pi/2)$,
\beaa
H^{\prime}(t) \LE \frac{2}{\sin ^{2} t} \left(-\cot t \cos ^{2} t-4 t - t\cot ^{2} t+4 t^{2} \cot t+2 t^{2} \cot ^{3} t\right)\\
\LT \frac{2}{\sin ^{2} t}\left(-\cot t \cos ^{2} t-t \cot ^{2} t+2 t^{2} \cot ^{3} t\right) \\
\EQ \frac{2\cot^{2}{t}}{\sin^{2}t} \left( -\frac{1}{2}\sin{2t} - t + 2t^{2}\cot{t} \right)\\
\LT 0,\eeaa
and
\[g^{\prime}(t)=\frac{t}{\sin ^{4} t}\left(3 \sin 2 t-2 t-4 t \cos ^{2} t\right)<0.\]
Moreover, $U$ has a unique zero
\be \lb{t0}t_0=t_{0,\im}:=U^{-1}(0)\in (0,\pi/2),\ee
and $g$ has also a unique zero
\be \lb{t1}t_1=t_{1,\im}:=g^{-1}(0)\in(0,\pi/2).\ee
Note that these two constants depend only on $(\im)$. Finally, it can be verified that $F$ has a unique zero
\be \lb{tau}\tau: = F^{-1}(0)\in(0,\pi/2).\ee

For the left end-point $t=\ta$, we have the following results.

\bb{lem}\lb{lem6} Suppose that $(\im)$ satisfies {\rm\x{im2}}. Then one has
\[\z\{ \ba{l}h_\imr(\ta_\imr) \geq 0  \iff r \geq P_\im\left( t_{0,\im}\right),\\
h_\imr(\ta_\imr) \le 0  \iff r \le P_\im\left( t_{0,\im}\right).
\ea\y.\]
Moreover, those equivalent inequalities take the equal sign at the same time. The functions $h_\imr,\ P_\im$ and constants $\ta_\imr,\ t_{0,\im}$ are as in {\rm\x{himr}, \x{Pim}}
and {\rm \x{thrr}, \x{t0}} respectively.
\end{lem}

\Proof Denote $h_\imr,\ P_\im,\ \ta_\imr,\ t_{0,\im}$ as $h,\ P,\ \ta,\ t_0$ respectively. Due to the relation \x{corr}, we modify the function $h(t)$ as
\[H(t)=\left(t-\cot t+\tpi \right)^{2}+ 4\left(t+\tpi \right) \cot t-t^{2} \left(2+\cot ^{2} t\right)^{2}.\]
By \x{corr}, \x{himr} and \x{Him},
\[H(\ta)=h(\ta).\]
By \x{Uim}, \x{Him} and \x{t0},
\[H(t_0)=0.\]
As $P(t)$ and $H(t)$ are strictly decreasing, we know that
$$r \ge P\left( t_{0}\right) \iff \ta=P^{-1}(r)\le t_0\iff h(\ta)=H(\ta) \ge H(t_0) = 0. $$
Similarly,
\[r\le P(t_0)\iff \ta\ge t_0\iff h(\ta) \le 0.\] Now the proof is finished. \qed

We also need some information on the derivative $h^{\prime}(t)$ at the end-point $t=\ta_\imr$. Notice that
\be \lb{h_diff}h^{\prime}(t)=2\left(2+\cot^{2}{t}\right) g(t).\ee

\bb{lem}\lb{lem4} Suppose that $(\im)$ satisfies {\rm\x{im2}}. Then
\[h_\imr^{\prime}(\ta_\imr)\geqslant 0
\iff r \geqslant P_\im\left( t_{1,\im}\right),\]
\[h_\imr^{\prime}(\ta_\imr)\leqslant 0
\iff r \leqslant P_\im\left( t_{1,\im}\right).\]
Moreover, those equivalent inequalities take the equal sign at the same time. The functions $g_\im,\ P_\im$ and constants $\ta_\imr,\ t_{1,\im}$ are as in {\rm\x{gim}, \x{Pim}}
and {\rm \x{thrr}, \x{t1}} respectively.
\end{lem}

\Proof We only prove the first equivalence relation and the other is similar. We omit $i,\ m$ in the subscripts of notations as in Lemma \ref{lem6}. As $P(t)$ is strictly decreasing, by \x{trtr} one knows that $r \ge P\left( t_{1}\right)$ if and only if $\ta=P^{-1}(r)\le t_1$. Then by the definition of $t_1$ and properties of $g$, it is also equivalent to $g(\ta) \geqslant 0$. \qed

For the above three important constants $\piim$, $P_\im(t_{0,\im})$ and $P_\im(t_{1,\im})$ relevant to the range of $r$ in Lemma \ref{lem5}, Lemma \ref{lem6} and  Lemma \ref{lem4}, we have the following observations.

\bb{lem} \lb{compare}  Suppose that $(\im)$ satisfies {\rm\x{im2}}. Then $P_\im\left( t_{1,\im}\right) < \piim$ and
 $P_\im\left( t_{1,\im}\right) < P_\im\left( t_{0,\im}\right)$. Moreover, with $\tau \in(0,\pi/2)$ in {\rm\x{tau},} we have
 \[P_\im\left( t_{0,\im}\right) -\piim \z\{\ba{ll} <0, \q t_{0,\im}\in(0,\tau),\\
 =0, \q t_{0,\im}=\tau,\\
 >0, \q t_{0,\im}\in(\tau,\pi/2).\ea \y.\]
 \end{lem}

\Proof Fix $i,\ m$ and denote $P_\im, \ t_{0,\im}, \ t_{1,\im}$ as $P,\ t_0,\ t_1$ respectively. From \x{t1} one has
$$\tpi= \cot{t_{1}} + \frac{t_{1}}{\sin{t_{1}^{2}}}- \frac{2t_{1}^{2} \cot{t_{1}}}{\sin^{2}{t_{1}}}.$$
By using this substitution and \x{Pim},
\[P\left( t_{1}\right) - \piim
= i^{2}\z( \tpi ^{2} + \z(4\cot{t_{1}} + \pi \y) \cdot \tpi  + 4 t_{1} \cot{t_{1}} - \frac{3}{4} \pi^{2} \y)=d(t_1),\]
where the function $d(t)$ is defined on $(0,\pi/2)$ as
\[d(t) :=
\left( \cot{t} + \frac{t}{\sin{t^{2}}} - \frac{2t^{2}\cot{t}}{\sin^{2}{t}} \right)
\left( 5\cot{t} + \frac{t}{\sin^{2}{t}} - \frac{2t^{2}\cot{t}}{\sin^{2}{t}} + \pi \right)
+ 4t\cot{t} - \frac{3}{4}\pi^{2}.\]
One can verify that $d(t)<0$ for $t\in(0,\pi/2)$. In particular, $d(t_1)<0.$

To prove $P_\im\left( t_{1,\im}\right) < P_\im\left( t_{0,\im}\right)$, since $P(t)$ is strictly decreasing, we only need to show that $t_1>t_0$. By \x{t0}, \x{t1} and the properties of functions $U$ and $g$, it is enough to hold that
\[\cot{t} + \frac{t}{\sin^{2}{t}} - \frac{2t^{2} \cot{t}}{\sin^{2}{t}} < t + t\cot^{2}{t} - \cot{t}, \q t\in(0,\pi/2),\]
which can be easily verified.

By using \x{Pim}, \x{t0} and replacing $\tpi$, one can get
\[P\left( t_0\right) - \piim= i^2 F(t_0),\]
where the function $F$ is as \x{Funct}. Compared with the zero $\tau$ of $F$ defined in \x{tau}, one has $F(t)<0$ on $(0,\tau)$ and $F(t)>0$ on $(\tau, \pi/2)$. Now the proof is finished. \qed

\bb{rem}\lb{cmp} By direct computations, it is easy to see that
\[F(0.5)<0,\qquad F(0.52)>0.\]Thus
$\tau\in (0.5, 0.52).$
Moreover, for $t_{0,\im}\in(0,\tau)\subset (0,0.52)$,
\[\tpi_\im < 0.52 + 0.52\cot^{2}{0.52} - \cot{0.52} < 0.36,\] by using {\rm\x{t0}} and the fact that $U$ is increasing.
Besides, $P_\im\left( t_{0,\im}\right) - \piim$ could be either negative or positive.
\end{rem}

Let us introduce a constant
\be \lb{r*}\ol{r}_\im: = P_\im\left( t_{0,\im}\right) >0.\ee

\bb{lem} \lb{lem9}
Let $m\in\N$, $1\le i \le m$ and $\ol{r}_\im$ be the constant given in {\rm\x{r*}}. If
\[i> \frac{m+1}{2}\andq r\ge \ol{r}_\im,\]
then
\[T_\imr =  \mct_\imr(r)=\bfw_\im^{-1}(r),\]
where the function $\bfw_\im$ is given in {\rm\x{WimT}}.
\end{lem}

\Proof We fix $i,\ m$ and omit them in the subscripts of all notations. We only need to prove that
\be \lb{111}\min_{\hr\in[0,r]}\mct(\hr) = \mct(\hr)|_{\hr=r}.\ee
From Lemma \ref{compare} we know that $r\ge P\left( t_{0}\right)> P \left(  t_1\right)$. By Lemma \ref{lem6} and Lemma \ref{lem4}, we get $h(\ta) \ge 0$ and $h^{\prime}(\ta) > 0$. According to \x{h_diff} and the properties of function $g$, one can deduce that $h$ is strictly increasing on $[\ta,  t_1]$ and strictly decreasing on $[ t_1, \pi/2]$.

We first consider the case
\[r\ge\piim.\]
It follows from Lemma \ref{lem5} that $h(\pi/2) \ge 0$ and thus $h(t)\ge 0$ on $[\ta, \pi/2]$, which together with Lemma \ref{signs} means that \x{111} holds.

Now we consider the case
$$\ol{r}_\im < \piim \andq r\in\z[\ol{r}_\im, \piim\y].$$
By Lemma \ref{lem5} one has $h(\pi/2) < 0$. It follows from Lemma \ref{signs} that $\mct (\hr)$, as a function of $\hr\in[0,r]$,  attains its minimum at $\hr=0$ or $\hr=r$. In fact, there exists $\hr^*\in(0,r)$ such that $\mct (\hr)$ is increasing on $[0,\hr^*]$ and decreasing on $[\hr^*,r]$.

Recall that $\mct(0)=\mct(\hr)|_{\hr=0}$ and $\mct(r) =\mct(\hr)|_{\hr=r}=\T$. In order to get \x{111}, we only need to prove that \[\mct(0) > \T= \bfw_\im^{-1}(r).\] As $\bfw(\d)$ is decreasing, it is equivalent to prove that
\be \lb{rc}r> \bfw(\mct(0)).\ee
Recall that
\[\tpi=\tpi_\im=\f{(m+1)\pi}{2i} -\f{\pi}{2},\quad\f{(m+1)\pi}{2i}= \tpi + \f{\pi}{2}.\]
Then
\[\mct(0)=\f{\pi}{\f{(m+1)\pi}{2i}+\sqrt{\z(\f{(m+1)\pi}{2i}\y)^2 + \f{r}{i^2}}}=\f{\pi}{\tpi +\f{\pi}{2} +\sqrt{\z(\tpi +\f{\pi}{2}\y)^2 + \f{r}{i^2}}},\]
\[\bfw_\im(\mct(0))=\f{4i^2\tpi}{1-\mct(0)} \cot\f{\tpi \mct(0)}{1-\mct(0)}
=K_{i,r}(\tpi) ,\]
where the function $K_{i,r}(u)$ is given as
\[K_{i,r}(u):= 4i^2 \z(u +\f{u \pi}{u-\f{\pi}{2} +\sqrt{\z(u+\f{\pi}{2}\y)^2 + \f{r}{i^2}}} \y) \cot\f{u\pi }{u-\f{\pi}{2} +\sqrt{\z(u+\f{\pi}{2}\y)^2 + \f{r}{i^2}}}.\]
Then the desired inequality \x{rc} is equivalent to
\be \lb{222}K_{i,r}(\tpi_\im) <r.\ee

According to Lemma \ref{compare} and Remark \ref{cmp}, one has $\tpi_\im \in (0, 0.36)$. For $u\in(0,0.36)$, we introduce another variable
\be \lb{vr}v =\frac{u\pi}{u-\f{\pi}{2} +\sqrt{\z(u+\f{\pi}{2}\y)^2 + \f{r}{i^2}}}.\ee
Note that, for $u\in (0,0.36)$, $v$ varies in
\[v\in
J_u:=\left[
\frac{u\pi}{u + \f{\pi}{2}},
\frac{u\pi}{u-\f{\pi}{2} +\sqrt{\z(u+\f{\pi}{2}\y)^2 + \f{r}{i^2}}}
\right]. \]
Thus the variables $(u,v)$ belong to the following domain
\[\Om=\Om_{i,r} := \z\{ (u,v): u \in \in(0,0.36), \ v\in J_u\y\}.\]
By using the variables $u$ and $v$ in \x{vr}, the function $K_{i,r}(u)$ is simplified to
\[K_{i,r}(u) = 4 i^2 (u+v) \cot v.\]
Solving $r$ from \x{vr}, we obtain
\[r= i^{2}\z(\frac{u^{2} \pi^{2}}{v^{2}} - 2u\pi - \frac{2u^{2}\pi}{v} + \frac{u\pi^{2}}{v}\y).\]
Thus inequality \x{222} can be deduced from the following one
\[4 i^2 (u+v) \cot v <i^{2}\z(\frac{u^{2} \pi^{2}}{v^{2}} - 2u\pi - \frac{2u^{2}\pi}{v} + \frac{u\pi^{2}}{v}\y).\]
Then \x{222} is guaranteed if
\be \lb{333}
f(u, v):= \frac{u^{2} \pi^{2}}{v^{2}} - 2u\pi - \frac{2u^{2}\pi}{v} + \frac{u\pi^{2}}{v} - 4(u + v) \cot{v} > 0,\quad (u,v)\in \Om.
\ee

To verify \x{333}, one has
\beaa
\frac{\partial{f(u,v)}}{\partial{v}} \EQ
-\frac{2u^{2}\pi^{2}}{v^{3}}
+ \frac{2u^{2}\pi}{v^{2}}
- \frac{u\pi^{2}}{v^{2}}
- 4\cot{v} + \frac{4(u + v)}{\sin^{2}{v}}\\
\LT -\frac{u\pi^{2}}{v^{2}} - 4\cot{v} + \frac{4(u + v)}{\sin^{2}{v}}\\
\EQ \frac{F(u, v)}{v^{2}\sin^{2}{v}},
\eeaa
where
$$F(u,v) = 4v^{2}(u + v) - 2v^{2}\sin{2v} - u\pi^{2}\sin^{2}{v}.$$
It is easy to verify that \[\frac{\partial F(u, v)}{\partial u} < 0\] for all $v\in(0,\pi/2)$.  From \x{Pim} and Remark \ref{cmp}, one has $ t_0\in(0,0.52)$ and
\[P ( t_0)=4 i^{2} \cot{t_{0}} \cdot
\left(t_{0} + \tpi  \right) > 4i^{2} \cot{0.52} \cdot 0.52 > 3.6i^{2}.\]
Then
\[\frac{u\pi}{u - \f{\pi}{2} + \sqrt{(u +\f{\pi}{2})^{2} + \frac{P ( t_0)}{i^{2}}}}<
\frac{u\pi}{u -\f{\pi}{2} + \sqrt{(u +\f{\pi}{2})^{2} + 3.6}},\]
and the domain
$$\tilde\Om=\z\{(u,v):u\in (0,0.36),\ v\in\z[
\frac{u\pi}{u +\f{\pi}{2}},\frac{u\pi}{u -\f{\pi}{2} + \sqrt{(u +\f{\pi}{2})^{2} + 3.6}}\y]\y\}\ \supset \Om.$$
It can be verified that
$$F\left(u, \frac{u \pi}{u -\f{\pi}{2}+ \sqrt{\left(u +\f{\pi}{2} \right)^{2} + 3.6}}\right) < 0,\quad  u \in (0, 0.36).$$
By considering the shape of the domain
$\tilde\Om$
and the properties of $F(u,v)$, one can get that $F(u,v)<0$ on $\tilde\Om$. Thus
\[\frac{\partial f(u, v)}{\partial v} < 0, \quad (u,v)\in \Om,\] and the left is to prove that inequality \x{222} holds when $r=P ( t_0)$.

Suppose now that $r=P(t_0)$. Since $P$ is strictly decreasing, it is equivalent to show that
\[t_0 <\frac{u\pi}{u -\f{\pi}{2} + \sqrt{\z(u +\f{\pi}{2} \y)^{2} + \frac{P ( t_0)}{i^{2}}}}.\]
By using $P ( t_0)=4 i^{2} \cot{t_{0}} \d\left(t_{0} + \tpi \right)$ from \x{Pim} and the substitution
$$u=\tpi =  t_0 +  t_0  \cot^{2}{ t_0} - \cot{ t_0} $$
from \x{t0}, one can verify that the inequality above holds when $ t_0\in(0,0.52)$. Now the proof is finished. \qed

Now we introduce the following two functions:
\[P_*(t)= i^2\z(t^{2}\left(2+\cot^{2}{t}\right)^{2}-\left(-\cot{t} + t + \tpi \right)^{2}\y),\]
and
\[V_\im(t)=\sqrt{(m+1)^2\pi^2+4P_*(t)}-(m+1)\pi-\frac{4it}{\pi}\left(t\cot^2t+\cot t+t-\tpi \y),\quad t\in(0,\pi/2).\]
By direct computations, we obtain that
\[V_\im'(t)=\frac{4i^2g(t)}{\pi\sqrt{(m+1)^2\pi^2+4P_*(t)}}\z(\frac{\sqrt{(m+1)^2\pi^2+4P_*(t)}}{i}-\pi(\cot^2t+2)\y)<0.\]
Note that
$$\sqrt{(m+1)^2\pi^2+4P_*(t)}/i-\pi(\cot^2t+2)<0,\quad t\in(0,\pi/2),$$
and the properties of the function $g$, we know that there exists $t_1\in(0,\frac{\pi}{2})$ such that $V_\im(t)$ is strictly decreasing on $(0,t_1)$ and strictly increasing on $(t_1, \pi/2)$.  As $\lim_{t\rightarrow0^+}V_\im(t)=\infty$ and $V_\im(\pi/2)=0$, there is a unique zero
\[\ul{t}_\im=\ul{t}: =V_\im^{-1}(0)\in (0,t_1) \subset (0, \pi/2).\]
We will use the notion
\be \lb{r_*}\ul{r}_\im=\ul{r}:=P_*(\ul{t}).\ee

\bb{lem} \lb{optvals}Assume that $(\im)$ satisfies {\rm\x{im2}}.  Let $\ol{r}_\im$ and $\ul{r}_\im$ be constants given by {\rm\x{r*}} and {\rm\x{r_*}}. Then for
$r\in(0, \ul{r}_\im],$
\[T_\imr = \frac{2i\pi}{(m+1)\pi+\sqrt{(m+1)^2\pi^2+4r}},\] while for $r\in(\ul{r}_\im,\ol{r}_\im),$ \[ T_\imr = \frac{2i^2 \si}{r} \left(\si\cot^2\si+\cot \si+\si-\tpi \y),\]
with $\si=\si_\imr$ being a constant given  as  {\rm\x{sigma}} below.
\end{lem}

\Proof For simplicity, we omit index $i,\ m$ in the subscripts of notations. From \x{PTr1} and Lemma \ref{signs}, $\mct_\imr$ can be seen as functions of $\hr\in[0,r]$ or $t\in[\ta,\pi/2]$ and we also simply denote them as $\mct_r(\hr)$ and $\mct_r(t)$. By Lemma \ref{compare}, $P(t_1)<P(t_0)=\ol{r}_\im$. As $P$ is strictly decreasing, one has $t_1>t_0$.

First we consider the case $r\le P(t_1)$. It follows from Lemma \ref{lem6} and Lemma \ref{lem4} that $h_r(\ta)<0$, $g(\ta)\le 0$. By \x{h_diff} and the fact that $g(t)$ is decreasing, one knows that $h_r(t)< 0$ on $[\ta, \pi/2]$. Then by Lemma \ref{TerP1}, Lemma \ref{signs} and \x{Thrr}, one has
\be\lb{tr0} T_\imr = \min_{\hr\in[0,r]}\mct_{r}(\hr) = \mct_{r}(\hr)|_{\hr=0} = \frac{2i\pi}{(m+1)\pi+\sqrt{(m+1)^2\pi^2+4r}}.\ee

Next we consider the case $P\left( t_1\right) <r< P\left( t_0\right)$. By Lemma \ref{lem6} and Lemma \ref{lem4}, we have $h_r(\ta)<0$ and $g (\ta) > 0$. One has $t_1>P^{-1}(r)=\ta$ due to \x{trtr} and the fact that $P$ is strictly decreasing. Thus according to \x{h_diff} and the properties of function $g (t)$, $h_r(t)$ is strictly increasing on $[\ta,  t_1]$ and strictly decreasing on $[ t_1, \pi/2]$, and of course \[\max_{t\in[\ta,\pi/2]}h_r(t) = h_r(t_1).\] Then from \x{himr} one knows that $h_r(t_1)$ is continuous and strictly increasing in $r$. We assert that $h_{P(t_1)}(t_1)<0$ and $h_{P\left( t_0\right)}(t_1)>0$. When $r=P(t_1)<P(t_0)$, one has $\ta=P^{-1}(r)=t_1$, then from Lemma \ref{lem6} there holds $h_r(t_1)=h_r(\ta)<0$. When $r=P(t_0)>P(t_1)$, one has $\ta=P^{-1}(r)=t_0<t_1$, then from Lemma \ref{lem6} it holds that $h_r(t_0)=h_r(\ta)=0$. Combining with the fact that $h_r(t)$ is strictly increasing on $[\ta, t_1]$, one can deduce that $h_r(t_1)>0$. Thus the assertion is true and there is some
\be \lb{r0}r_0:=r_{0,\im}\in\z(P\left( t_1\right), P\left( t_0\right)\y)\ee
such that $h_{r_0}(t_1)=0.$
If $P\left( t_1\right)<r\le r_0$, then $h_r(t)\le 0$ on $[\ta, \pi/2]$ and \x{tr0} holds obviously.

Therefore, in the following we assume that $r_0<r<P\left( t_0\right)$. Then $h_r(t_1)>0$ and $h_r(t)$ has a unique zero $\si$ in $(\ta,t_1)$ as well as in $(0,t_1)$. We denote the function $h_r(t)$ defined in interval $(0,t_1]$ as $h|_{(0, t_1]}(t)$. Let
\be \lb{sigma}\si=\si_\imr:=h|_{(0,t_1]}^{-1}(0)\in (\ta,t_1),\ee
which is a constant decided by the values of $i,\ m,\ r$. From Lemma \ref{signs}, one can easily get that $\mct_{r}(t)$ is strictly decreasing on $[\ta, \si]$ and then decrease monotonically on $[\si,\pi/2]$, which is relevant to the sign of $h_r(\pi/2)$. By Lemma \ref{TerP1} and Lemma \ref{signs}, we know that
\[T_\imr = \min_{\hr\in[0,r]}\mct_{r}(\hr) = \min_{t\in[\ta,\pi/2]}\mct_{r}(t) = \min\Big\{\mct_{r}(t)|_{t=\si},\mct_{r}(t)|_{t=\pi/2}\Big\}.\]
By considering \x{tr0} and function $t(\hr)\ \z(=t_\im(\hr)\y)$ in Lemma \ref{signs}, one has
\be \lb{tpi/2}
\mct_{r}(t)|_{t=\pi/2} = \mct_{r}(\hr)|_{\hr=0} = \frac{2i\pi}{(m+1)\pi+\sqrt{(m+1)^2\pi^2+4r}}.\ee
As $h_r(\si)=0$, from \x{himr} we have
\be \lb{rsi}r=P_*(\si).\ee

 Since $P_*(t)=-i^2h_r(t)+r$, we can deduce the monotonicity of $P_*$ and get that \[\min_{t\in(0,\pi/2)}P_*(t) =P_*(t_1).\] One can easily prove $P_*(t_1)>0$ by using the substitution
$$\tpi= \cot{t_{1}} + \frac{t_{1}}{\sin{t_{1}^{2}}}- \frac{2t_{1}^{2} \cot{t_{1}}}{\sin^{2}{t_{1}}}$$
according to \x{t1}. Then $P_*(t)>0$ for all $t\in(0,\pi/2)$. By \x{trr2} and \x{rsi},
\be \lb{tsi}\mct_{r}(t)|_{t=\si}=\frac{2i^2 \si}{r} \left(\si\cot^2\si+\cot \si+\si-\tpi \y).\ee

By using \x{tpi/2}, \x{rsi} and \x{tsi}, there holds
\be \lb{dif}\mct_{r}(t)|_{t=\pi/2}-\mct_{r}(t)|_{t=\si}=\frac{i\pi}{2r}V_\im(\si).\ee
By \x{dif}, $\mct_{r}(t)|_{t=\pi/2}>\mct_{r}(t)|_{t=\si}$ when $\si<\ul{t}$ and $\mct_{r}(t)|_{t=\pi/2}<\mct_{r}(t)|_{t=\si}$ when $\si>\ul{t}$. Besides, they take the equal sign at the same time.  As $\si, \ \ul{t}\in(0,t_1)$ and $P_*(t)$ is strictly decreasing in $t\in(0,t_1)$, we know that $\ul{t}-\si$ has the same sign with $P_*(\si)-P_*(\ul{t})=r-\ul{r}$.

Finally, we show that $\ul{r}\in\z(r_0, P\left( t_0\right)\y)\subset\z(P\left( t_1\right), P\left( t_0\right) \y)$, where $r_0$ is as in \x{r0}. Otherwise if $\ul{r}\le r_0$, from the assumption $r_0<r<P\left( t_0\right)$ one has $r>\ul{r}$. Then there hold $\ul{t}>\si$ and
\[\mct_{r}(t)|_{t=\pi/2}>\mct_{r}(t)|_{t=\si}.\] From the discussion above and the fact that $\mct_r$ and $\si$ are continuous in $r$, one can deduce that
\[\mct_{r_0+\epsilon}(t) |_{t=\pi/2}<\mct_{r_0+\epsilon}(t) |_{t=\si}\] when $\epsilon>0$ is small enough. This is a contradiction. In a similar way, we can show that
$\ul{r}\ge P\left( t_0\right)$ is also impossible. \qed

Combined Lemma \ref{ile}, Lemma \ref{lem9} and Lemma \ref{optvals}, we have proved the following result.

\bb{thm} \lb{optpts} For $r\in(0,\oo)$, $m\in\N,~1\le i \le m,$ let
 $\ol{r}_\im$ and $\ul{r}_\im$  be constants given by  {\rm\x{r*}}  and {\rm\x{r_*}} respectively.
We have the following results.
\begin{itemize}
\item [{\rm (i)}]If $i\le \frac{m+1}{2}$, then $$T_\imr=\mct_\imr(\hr)|_{\hr=r} =\bfw_\im^{-1}(r),$$
where $\bfw_\im$ is given by {\rm\x{WimT}}.
    \item[{\rm (ii)}] If $i>\frac{m+1}{2}$ and $r\in [\ol{r}_\im,\oo)$, then $$T_\imr=\mct_\imr(\hr)|_{\hr=r}=\bfw_\im^{-1}(r).$$
\item[{\rm (iii)}] If $i>\frac{m+1}{2}$ and $r\in(\ul{r}_\im,\ol{r}_\im)$, then $$T_\imr=\mct_\imr(\hr)|_{\hr=\hr(\si)}= \frac{2i^2 \si}{r} \left(\si\cot^2\si+\cot \si+\si-\frac{m-i+1}{2 i} \pi\y),$$
where $\si$ is given by {\rm\x{rsi}} and
    \be \lb{hrsi}
    \hr(\si)=2i^2 \cot \si \cdot\z(\si\cot^2\si-\cot \si+3\si+\tpi \y).
    \ee
 \item [{\rm (iv)}]If $i>\frac{m+1}{2}$ and $r=\ul{r}_\im$, then $$T_\imr=\mct_\imr(\hr)|_{\hr=0}=\mct_\imr(\hr)|_{\hr=\hr(\si)} = \frac{2i\pi}{(m+1)\pi+\sqrt{(m+1)^2\pi^2+4r}}.$$
\item[{\rm (v)}] If $i>\frac{m+1}{2}$ and $r\in(0, \ul{r}_\im)$, then $$T_\imr=\mct_\imr(\hr)|_{\hr=0}= \frac{2i\pi}{(m+1)\pi+\sqrt{(m+1)^2\pi^2+4r}}.
    $$
     \end{itemize}
\end{thm}

Up to now, we have worked out the solutions for the minimization problem \x{Timr}. Next we show that such solutions are just those for
our minimization problem \x{tmin}, based on the strong continuity stated as in Lemma \ref{CC}.

\subsection{Solutions for the minimization problem \x{tmin}}
For $m\in\N$, $1\le i \le m$, $\hr,\ \kr\ge 0$ and $T\in(0,1)$. Let us introduce the measure
$\gamma_{i,m,\hr,\kr,T}$ as
\bea \lb{gamma}
\gamma_{i,m,\hr,\kr,T}(x)\EQ\z\{\ba{l} \mu_{i-1,\hr}(\frac{x}{T}), \qq\qq\qq\qq\q  x\in[0,T],\\
-\hr+\frac{1}{1-T}\int_T^x\eta_{m-i,\kr}(\frac{t-T}{1-T})\dt, \q  x\in(T,1],
\ea \y.\eea
where  $\mu_{i-1,\hr}$ and $\eta_{ m-i, \kr}$ are given as in \x{mumr} and  in \x{etamr}, respectively.
Note that the constraint of the measure $\gamma_{i,m,\hr,\kr,T}(x)$ on $[T,1]$ is actually the scaling measure $\mu_{\eta_{ m-i, \kr}}(\frac{x-T}{1-T})$ obtained by \x{muq}.

\bb{lem} \lb{criticalmeasurelem} For $r\in(0,\oo)$, $m\in\N,~1\le i \le m$ and define the measure
$\gamma_{i,m,r}(x)$ as
\[\gamma_{i,m,r}= \z\{\ba{ll} \gamma_{i,m,r,0,T_\imr}, &  i\le \frac{m+1}{2},
\\ \gamma_{i,m,r,0,T_\imr}, & i>\frac{m+1}{2}, r\in [\ol{r}_\im,\oo),
\\  \gamma_{i,m,\hr(\si),r-\hr(\si),T_\imr},  & i>\frac{m+1}{2}, r\in(\ul{r}_\im,\ol{r}_\im),
\\\gamma_{i,m,0,r,T_\imr} \mbox{ or } \gamma_{i,m,\hr(\si),r-\hr(\si),T_\imr}, & i>\frac{m+1}{2} , r =  \ul{r}_\im,
\\ \gamma_{i,m,0,r,T_\imr}, & i>\frac{m+1}{2}, r\in(0, \ul{r}_\im). \ea \y.\]
 Then $ T_\imr= T_\im\z(\gamma_{i,m,r}\y).$
\end{lem}

\Proof By \x{L1rBT},\x{M1rBT} and \x{P2}, we have
 \[\la_{m}\z(\gamma_{i,m,r}(x)|[0,1] \y) = \la_{i-1}\z(\gamma_{i,m,r}(x)|[0,T_\imr] \y) =  \la_{m-i}\z(\gamma_{i,m,r}(x)|[T_\imr,1]\y).\]
Hence, $T_\imr= T_\im\z(\gamma_{i,m,r}\y).$
\qed

\bb{lem} \lb{L1s} For $r\in(0,\oo)$, $m\in\N$,  $1\le i \le m$ and $T_\imr$ as above. Then
\[\mimr = T_\imr.\]
\end{lem}

\Proof For any potential $q\in \sor$, we have constructed two parameters $\hr, \ \kr$  as in \x{r1r2}-\x{ine} such that the node $T=T_\im(q)$ satisfies the inequality \x{ine}. By using the function $f(t,\hr)$ defined in \x{ftrfun},
\[f(T_\im(q),\hr) \le 0= f(\mct_\imr(\hr),\hr).\]
Since $f(t,\hr)$ is decreasing in $t$, we obtain
\[T_\im(q) \ge \mct_\imr(\hr).\]
Hence $T_\im(q)\ge T_\imr$ for each $q\in S_r$. Therefore $\mimr=\inf_{q\in \sor} T_\im(q)\ge T_\imr$.\qed

On the other hand, by Lemma \ref{wstarlemma}, there exists a sequence $\{\mu_n\} \subset \Mo  \cap \mcc^\oo$
such that $ q_n= \mu_n'\in B_r$  and
\[\mu_{n} \to \gamma_{i,m,r}  {\rm ~in}~ (\Mo ,w^*),\quad n \to +\oo.\]
By the  continuity of $T_\im(\mu)$ in the measure $\mu$, we obtain that
 \[  T_\im(\gamma_{i,m,r}) = \lim_{n \to \infty} T_\im( \mu_n) = \lim_{n \to\infty} T_\im(q_n)
 \ge \lim_{n \to \infty} \mimr = \mimr.\]
Thus $  \mimr \le  T_\imr$ by Lemma \ref{criticalmeasurelem}. As a consequence, the result has been proved.
\qed

It follows from Theorem \ref{optpts} and  Lemma \ref{criticalmeasurelem}, Lemma \ref{L1s} that the main results of this paper
can be stated as follows.

\bb{thm}\lb{main}
For $r\in(0,\oo)$, $m\in\N$ and $1\le i \le m$. Let $\ol{r}_\im,~\ul{r}_\im,~\si=\si_\imr,~\hr(\si)$  be constants given as  {\rm\x{r*}, \x{r_*}, \x{sigma}, \x{hrsi}} respectively,
and  $\gamma_{i,m,\hr,\kr,T}(x)$ be the measure as {\rm\x{gamma}.} Then we have the following complete characterizations for $\mimr$:
\begin{itemize}
\item[{(i)}] If  $i\leq \frac{m+1}{2},$ then
\[\mimr = \bfw_\im^{-1}(r),\]
where $\bfw_\im$ is given by {\rm\x{WimT}}. Moreover, $\mimr$ can be attained by the measure $\gamma_{i,m,r,0,T_\imr}(x).$
 \item[{(ii)}] If  $i>\frac{m+1}{2}$ and $r\in [\ol{r}_\im,\infty),$ then
\[\mimr = \bfw_\im^{-1}(r).\] Moreover, $\mimr$ can be attained by the measure $\gamma_{i,m,r,0,T_\imr}(x)$.
    \item[{(iii)}] If $i>\frac{m+1}{2}$ and $r\in(\ul{r}_\im,\ol{r}_\im)$, then
\[\mimr = \frac{2i^2 \si}{r} \left(\si\cot^2\si+\cot \si+\si-\frac{m-i+1}{2 i} \pi\y).\]
Moreover,  $\mimr$ can be attained by the measure $\gamma_{i,m,\hr(\si),r-\hr(\si),T_\imr}(x)$.
\item[{(iv)}] If $i>\frac{m+1}{2}$ and $r\in(0, \ul{r}_\im]$, then
\[\mimr = \frac{2i\pi}{(m+1)\pi+\sqrt{(m+1)^2\pi^2+4r}}.\]
Moreover, $\mimr$ can be attained by the potential
\[q_\imr(x)=\z\{\ba{l}
0, \qq\qq\qq\qq\qq\q x\in[0,T_\imr],\\
\frac{1}{1-T_\imr}\eta_{m-i,r}\z(\frac{x-T_\imr}{1-T_\imr}\y), ~ x\in[T_\imr,1],\ea \y.\]
with $\eta_{m,r}(x)$ being the potential as in {\rm\x{etamr}.} Moreover, $\mimr$ can be also attained by the measure $\gamma_{i,m,\hr(\si),r-\hr(\si),T_\imr}(x)$ when $r=\ul{r}_\im$.
\end{itemize}\end{thm}

\bb{rem}{\rm The approach used in this paper can be applied to the following more general spectral
problem
\be\lb{gee}y''=q(x) y +\lambda m(x)y,\ee
where $q$ is nonnegative continuous function and $m$ may change sign.
See \cite{cg} for the distribution of eigenvalues for problem \x{gee} and \cite{cmz-adv} for the sharp bounds
of its eigenvalues. When $q(t) \equiv 1/4$, it is a spectral equation related to the famous Camassa-Holm equation.
See \cite{c-phd} for discussions on such a spectral problem. }
\end{rem}

\vspace{3mm}

\noindent {\bf Conflict of interest} The authors have no conflict of interest.

\bibliographystyle{amsplain}

\end{document}